\documentclass{llncs}
\usepackage[german,english]{babel}
\usepackage{times,theorem}
\usepackage{latexsym,color,graphics}

\usepackage{epsfig}
\usepackage{amssymb}
\usepackage{amsmath}
\usepackage{amsfonts}
\usepackage{xspace}
\usepackage{graphicx}
\begingroup
\catcode`\~=11
\gdef\urltilde{\lower 0.6ex\hbox{~}}
\endgroup

\newcommand{\A}{\mathcal{A}}

 \renewcommand{\L}{\mathcal{L}}

\title{General Categorial Geometry and Algebraic Topology}
\author{Zoran Majki\'c}
\authorrunning{Zoran Majki\'c}
\institute{ISRST, Tallahassee, FL, USA\\
\email{majk.1234@yahoo.com}}

\newtheorem{coro}{Corollary}

\begin{document}

\maketitle

\begin{abstract}
In Categorial Topology, given a category (as a "geometric object")  we can consider its properties preserved under continuous action (a "deformation") of a comma-propagation operation.
However, the Metacategory space, valid for all  categories,  cannot be defined by using well-know Grothendeick's approach with discrete ringed spaces.
So, we can consider any category $\textbf{C}$ as an abstract geometric object,
 that is, a discrete space where the points  are the objects of this category and morphisms between objects as the oriented paths.
  \\Based on this approach, we  define the Cat-arrows space $V$ valid for all categories with commutative (and associative) partial addition operation $\oplus$ for the vectors, based on partial operation of categorial composition of morphisms, their inner and outer products in 3D Cat-arrows space, like in 3D Clifford algebra. We provide a general definition of the norm ("weight") of the vectors in $V$, cumulative under the composition of the category morphisms. This norm assigned to the morphisms, non existing in metacategory theory, enriches each individual category with new semantical interpretation of the morphisms: they are "weighted" with a real number. So, this general transformation of metacategory theory into 3D Cat-arrows space of vectors is not only another example of applicability of 3D Clifford algebra, but the new semantic (of their "weight") and topological enrichment of the Metacategory theory with its abstract undefined primitive concept of morphisms.
\end{abstract}

\section{Introduction to Metacategory and its Geometric Space\label{sec:SymMet}}
 Category theory and geometry are deeply intertwined, with category theory providing a general framework for studying mathematical structures like geometric spaces by focusing on relationships (morphisms) between objects rather than their internal details. This approach is evident in areas like algebraic topology, where functors translate between geometric spaces and algebraic objects like groups.

 It is well known that, in a given category $\textbf{C}$, its morphisms (arrows) in $Mor_\textbf{C}$ can be composed by the non-commutative \emph{partial} operation $\circ:Mor_\textbf{C}\times Mor_\textbf{C} \rightarrow Mor_\textbf{C}$  to obtain more complex (composed) morphisms, so that in the associative non-commutative partial algebra $\A_M = (Mor_\textbf{C},\circ)$, for each object $a \in Ob_\textbf{C}$,  the identity morphism $id_a \in Mor_\textbf{C}$ is the "representation" of this object. From this point of view, the properties of the category can be represented by using \emph{only} its morphisms.
Such compositional property for the objects and the possibility to represent the category \emph{only} by using the objects of a category generally does not hold.

In what follows I will use the notion \emph{metacategory}  as it was defined in \cite{McLn71}, that is, by a metagraph consisting of objects $a,b,c,...$, morphisms (or morphisms) $f,g,h...$ and two operations, \emph{domain} which assigns to each morphism $f$ an object $a = dom(f)$, and \emph{codomain} which assigns to each morphism $f$ an object $b = cod(f)$.

Hence, a \emph{metacategory}  is a metagraph with two additional operations:
\begin{enumerate}
  \item \emph{Identity} (or $id$) which assigns to each object $a$ an morphism $id(a)$ with
 $$dom(id(a)) = cod(id(a)) = a$$
 denoted by $id_a:a\rightarrow a$ 
 %
\item \emph{Composition} which assigns to each pair of morphisms $(f,g)$ with $dom(g) = cod(f)$ an morphism $g\circ f$ called \emph{composite}.
\end{enumerate}
 These operations in a metacategory are subject to the following axioms:

1. \emph{Associativity} of composition, $k\circ(g\circ f) = (k\circ g)\circ f$, for each well defined composition of three morphisms, that is, an axiom in FOL:
\begin{equation}\label{eq:Assoc}
(\forall f)(\forall g)(\forall k)( ((cod(f) = dom(g))\wedge (cod(g) = dom(k)))\Rightarrow (k\circ(g\circ f) = (k\circ g)\circ f))
\end{equation}

2. \emph{Unit law}: for each morphism $f:a\rightarrow b$ and $g:b\rightarrow c$ composition with the identity morphism $id_b$ gives $id_b\circ f = f$ and $g\circ id_b = g$, that is, an axiom in FOL:
\begin{equation}\label{eq:UnitLow}
(\forall f)(\forall g)( (cod(f) = dom(g))\Rightarrow ((id(cod(f))\circ f = f)\wedge (g\circ id(dom(g)) = g)))
\end{equation}
That is, this notion of a metacategory is used for a category described directly by means of axioms without using the theory of sets.
A category (as distinguished from a metacategory) is to be any interpretation of these axioms within set theory. Commutative diagrams in a given category $\textbf{C}$ are expressed algebraically by a number of equations, such that
 \begin{itemize}
   \item  A \emph{normalized} commutative diagram in $\textbf{C}$ is obtained by elimination from it all identity morphisms.
\item An \emph{atomic} commutative diagram in $\textbf{C}$ is  expressed algebraically by a single equation. Thus, it can be expressed by a single morphism in the arrow category $\textbf{C}\downarrow\textbf{C}$.
 \end{itemize}
\textbf{Remark}: Note that also in definition of a metacategory, the unique operation for the objects is the \emph{Identity} operation which simply tells to us the existence of the identity morphism for each object of a given category. So, this operation does not tell anything \emph{about the objects} but from the fact that it is an operation on objects, we specified only that for a given object \emph{there is a unique} identity morphism of this object. \emph{Identity} operator is necessary only for definition of the Unit law in the axiom (\ref{eq:UnitLow}).

 Unit law in the case when $f = g =id_b:b\rightarrow b$ gives $id_b\circ id_b = id_b$ for each identity morphism, so that for any  object $b$  the set of all morphisms from $b$ into $b$ is a monoid with its unit equal to this identity morphism $id_b$. Again, the identity morphism does not tell us nothing about the  object from which is derived by operation \emph{Identity}, but tells us that it can be eliminated from the composition of the morphisms (for example, in order to obtain the \emph{normalized} commutative diagrams).

 Any \emph{category} is a specific (particular) collection of the objects and the morphisms (arrows).
A category is \emph{connected} if it is possible to go from any object to any other object of the category along a path of "composable" forward or backward morphisms. So, every category is a union of disjoint connected
subcategories in a unique way.
\\$\square$\\
A topos \cite{BaWe85,LaSc86}is a special kind of category defined by axioms saying roughly that certain constructions one can make with sets can be done in the category. In that sense, a topos is a generalized set theory. However, it originated with Grothendieck and Giraud as an abstraction of the properties of the category of sheaves of sets on a topological space. So, there exists an intimate connection between toposes and algebraic geometry \cite{Weil46}.

Later, Lawvere and Tierney introduced a more general idea which they called “elementary topos” (because their axioms were first order and involved no quantification over sets), and they and other mathematicians developed the idea that a theory in the sense of \emph{mathematical logic} can be regarded as a topos, perhaps after a process of completion.  Lawvere’s original insight was that a
mathematical theory (corresponding roughly to the definition of a class of mathematical objects) could be usefully regarded as a category with structure of a certain kind, and a model of that theory  as a set-valued functor from that category
which preserves the structure. The structures involved are more or less elaborate, depending on the kind of objects involved. The most elaborate of these use categories which have all the structure of a topos.

As far as I know, there is the following approach in the category theory to consider its geometric properties, which shortly is provided in \cite{Cart01} (Section 5. Points and Representation) where representations corresponds to the morphisms of a given category and each object of this category is a space composed by the elements called the "points" as well). In this approach, the "points" are just particular morphisms from the object \textbf{1}:

"\emph{We introduce a space consisting of a single point \textbf{1}. For each point a of a space X there exists a unique mapping, or representation, of the space \textbf{1} into the space X mapping the one point of \textbf{1} to the point a of X. Category theory is the mathematical expression of the idea of representation (or transformation). We have a class of objects (or spaces) and transformations $f$ of an object X into an object Y, with the possibility of composing these transformations. The “idea” of category theory is to consider only spaces and transformations rather than points. However, in the majority of categories there exists an object such as \textbf{1} characterized by the fact that there is exactly one transformation from X into \textbf{1}, whatever the object X. One can thus call any transformation of \textbf{1} into X a point of X, but it may very well happen that an object has no points in this sense.}"

This approach is used, for example in \cite{Gold79}, to model a propositional logic (with the Boolean algebras representation). That is, it is  valid for a number of applications of Category Theory and definition of particular categories as models of different algebraic problems as, for example,  the \textbf{Set} category as topos of sets (in which each object is a set of elements (the "points") with  logical object $\Omega$ whose elements are two truth values (true/false), and the morphisms are the functions between the sets),  or the Cartesian Closed Categories (CCC) as a categorial model for the typed lambda calculus \cite{AsLo91}, or Grothedieck category $\mathfrak{T}(X)$ of sheaves over a space $X$ considered as a topos where the truth values correspond to the open sets of $X$.
A Grothendieck's topology is a structure on a category $\textbf{C}$ that makes the objects of $\textbf{C}$ act like the open sets of a topological space.
A category together with a choice of Grothendieck's topology is called a \emph{site}.

Grothendieck's topologies axiomatize the notion of an open cover.  Using the notion of covering provided by a Grothendieck topology, it becomes possible to define sheaves on a category. This was first done in algebraic geometry. Grothendieck's insight was that the definition of a \emph{sheaf} (introduced by J. Leray in the field of algebraic topology \cite{Houz94}) depends only on the open sets of a topological space, not on the individual points.
Grothendieck's refoundation of algebraic geometry and the panoply of new notions
that he introduced (such as schemes, topoi, motives, and stacks) constitute one of the most important paradigm shifts that occurred in mathematics during the XXth century.
\\
"\emph{Grothendieck did not derive his inspiration from physics
and its mathematical problems. Not that his mind was incapable of grasping this
area|he had thought about it secretly before 1967 - but the moral principles that
he adhered to relegate physics to the outer darkness, especially after Hiroshima. It is surprising that some of Grothendieck's most fertile ideas regarding the nature of space and symmetries have become naturally wed to the new directions in modern physics.}" p.389 in \cite{Cart01}.

Moreover, Grothendieck’s approach is strongly based on (and significantly contributed
to develop) category theory and groupoid theory, thereby having a seminal impact on the “foundations of mathematics”.
\\
However, in our case, by considering what can be the general Metacategory space, i.e., the space of  \emph{each} given category $\textbf{D}$, with its family of objects $X = Ob_\textbf{D}$ and family of arrows (morphisms) $Mor_\textbf{D}$, by considering the  representation of categorial properties  by diagrams, expressed by oriented graphs which nodes in the plain space (real \emph{physics surface}) represented by positions (points) of the objects win $X$ and its oriented paths (the arcs of a graph) represented by the morphisms in $Mor_\textbf{D}$.

So, by using the Grothendieck's approach, where the space of category $\textbf{D}$ is a discrete  space $X$, for which any subset $U\subseteq X$ is an open space, we are able to define a new category $\textbf{C}$ whose objects in $Ob_\textbf{C}$ are the open sets of $X$  (all subsets of the set of objects $X$ of our category $\textbf{D}$, by considering that also each singleton in this discrete \emph{topological} space of the category $\textbf{C}$ is an open space as well), and the morphisms in  $Mor_\textbf{C}$ are the inclusion maps (functions) of on open set into another, just as in the context of \emph{schemes} (the collection of objects of $\textbf{C}$ is a topological space together with \emph{commutative rings} for all of its open (sub)sets (the addition is exclusive union with empty set as additive unity, while the multiplication is the set intersection with multiplicative unity $X$), so it is a scheme, i.e., \emph{ringed space}).

So, we are able to define the \emph{presheaf of} $X$ as a contravariant functor from $\textbf{C}$ to the category of sets $\textbf{Set}$, i.e., $F:\textbf{C}^{OP}\rightarrow \textbf{Set}$, where $\textbf{C}^{OP}$ has the same objects as $\textbf{C}$ but with opposite morphisms.
Open sets $\{U_i\}_{i \in I}$ \emph{covers} $U~$ iff $~\bigcup_{i \in I}  U_i = U$.
Functor $F$ maps any inclusion morphism $g:U_i\hookrightarrow U$  in $\textbf{C}$ into a restriction map  $F(g^{OP}): F(U) \rightarrow F(U_i)$  in $\textbf{Set}$, and for any element (section) $x\in F(U)$ we denote by $x|U_i$ the element  $F(g^{OP})(x)\in F(U_i)$.

  A \emph{sheaf} is a presheaf $F$ that satisfies the gluing axiom: if for such covering open sets $U_i$ each section $x_i \in F(U_i)$ is given for each $i\in I$ in such a way that $x_i|U_i\bigcap U_j = x_j|U_i\bigcap U_j$ for all $i,j \in I$, then there is a \emph{unique} $x\in F(U)$ such that $x|U_i = x_i$.

 We can replace each $U_i$
 with a \emph{family} of open subsets; in this example, $U_i$
is replaced by a family of  open immersions
$g:V_j \hookrightarrow U_i$,
\begin{equation}\label{eq:sheaf}
U_i ~~~~\mapsto ~~~~S =\{g:V_j \hookrightarrow U_i\}_{j \in I}
\end{equation}
Such a collection $S$ is called a \emph{sieve} if  for any morphism $g:V_j\rightarrow U_i$ is in $S$ and  $f:c'\rightarrow V_j$ is any other morphism in $\textbf{C}$, then $g\circ f$ is in $S$.
Consequently, sieves are similar to right ideals in ring theory or filters in order theory. If $S$  in (\ref{eq:sheaf}) is sieve such that  $\bigcup_{j \in I}\{V_j\} = U_i$, then $S$ is a \emph{covering} sieve.
Sieves were introduced by Giraud \cite{Gira64} in order to reformulate the notion of a Grothendieck's topology: the notion of a collection of open subsets of $U$ stable under inclusion is replaced by the notion of a sieve.

If $c$ is any given object in $\textbf{C}$, a \textbf{sieve} on $c$ is a subfunctor $F \in \textbf{Set}^{\textbf{C}^{OP}}$ of the functor $Hom(\_,c)\in \textbf{Set}^{\textbf{C}^{OP}}$  (this is the Yoneda embedding applied to $c$, as shown in the commutative diagram provided by Figure 1 in Appendix (Section 6)).

From the fact that that $\textbf{C}$ is locally small (at maximum we can have one inclusion map between any two objects in $\textbf{C}$), we can apply the Yoneda lemma as well, so that for fixed this particularly relevant object $c$ in $\textbf{C}$ (i.e., in $\textbf{C}^{OP}$), the natural transformations $Nat(Hom(\_,c), F)$ is one-to-one correspondence with the elements of $F(c)$ (an “element” of $F(c)$ is therefore a natural transformation into sieve $F$ on $c$, that is there is Yoneda bijection
$$\varphi:Nat(Hom(\_,c), F) \backsimeq F(c)$$
and the  diagram in $\textbf{Set}$, provided by Figure 1 in Appendix (Section 6), commutes for any monomorphism $g:c'\hookrightarrow c$ in $\textbf{C}$ with the function $Hom(g^{OP},c) = \_\circ g$ such that $Hom(g^{OP},c)(id_c) = id_c\circ g = g$, for a  natural transformation $\phi \in Nat(Hom(\_,c), F)$,
where, from the fact that $F$ is a sieve (a subfunctor of $Hom(\_,c)$), $F(c)\subseteq Hom(c,c)$,  $F(c')\subseteq Hom(c',c)$ and $F(g^{OP})$ is a restriction of $Hom(g^{OP},c)$ to $F(c)$.  Then $c'$ will be considered "selected" by $F$ if and only if $F(c')$ is nonempty.

The most common operation on a sieve is pullback. Pulling back a sieve $F$ on $c$ (a subfunctor of $Hom(\_,c)$ in $\textbf{Set}^{\textbf{C}^{OP}}$)  by a morphism $g:c'\hookrightarrow c$ in $\textbf{C}$ gives a new sieve $g^*F$ on $c'$. This new sieve is called "pullback of $F$ along $g$", defined as the fiber product in category of functors $\textbf{Set}^{\textbf{C}^{OP}}$, $F\times_{Hom(\_,c)}Hom(\_,c')$, together with natural transformation embedding (projection) $p_2$ in $Hom(\_,c')$,
natural transformation embedding $p_1$ in $F$, natural transformation $in$ representing that $F$ is a subfunctor of $Hom(\_,c)$ and natural transformation $Hom(\_,g)$, as shown by Figure 2 in Appendix (Section 6). 
That is, with this we expressed a Grothendieck's topology (a collection of \emph{covering} sieves) by the following axioms:
\begin{itemize}
  \item First axiom of \emph{base change}: If $F$ is a covering sieve on $c$  and $g: c' \hookrightarrow c$ is a morphism, then the pullback $g^*F$ is a covering sieve on $c'$ (i.e., if $\{U_i\rightarrow c\}_{i\in I}$ covers $c$ and $g: c' \hookrightarrow c$ then $\{U_i \times_c c'\rightarrow c'\}_{i\in I}$ covers $c'$ (i.e., the pullback $g^*F$ is a covering sieve on $c'$) where fibre products $\times_c$ in $\textbf{C}$ are intersections.
  \item Second axiom is of \emph{local character} (if $\{U_i\rightarrow c\}_{i\in I}$ covers $c$  and $\{V_{ij}\rightarrow U_i\}_{j\in J}$ covers $U_i$ for each $i$, then the collection $\{V_{ij}\rightarrow U_i\rightarrow c \}$ for all $i$ and $j$ should cover $c$).
  \item The \emph{identity} axiom: $Hom(\_, c)$ is a covering sieve on $c$ for any object $c$ in $\textbf{C}$ (any set is covered by itself via the identity map).
\end{itemize}
  By introducing the category theoretic notion of \emph{functor of solutions} (or, in the dual geometric version, functor of points \cite{Vacu21}), Grothendieck
takes into account the whole family of “figures” at once.

That is, in both cases, if the functor $F$ is a shief or sieve, it is the representation of the objects of $\textbf{C}$  in the category of sets $\textbf{Set}$ (Grothendieck's representation theorem states
“that every ring is isomorphic to the ring of global sections of a sheaf of local rings”).\\
\textbf{Remark}: \emph{However}, no one of such geometric representations can  represent the space of our given category $\textbf{D}$: the derived from it category $\textbf{C}$ has as morphisms the inclusion maps between open sets considered as spaces composed by the points (in the case when the functor $F$ is a sheaf, each objects  of $\textbf{C}$, see in (\ref{eq:sheaf}), is a subset of the objects of $\textbf{D}$, while if $F$ is a sieve, each objects of $\textbf{C}$ are replaced by a collection of inclusion maps (right hand side in (\ref{eq:sheaf})).
\\$\square$\\
Thus, in both  Grothendieck's topology approaches \cite{GrDi60}, the geometric space represented by them do not represent the morphisms of the category $\textbf{D}$ (considered by me as the paths in the $\textbf{D}$ category discrete topological space, and do not represent a point of such space as a single point of it). This is the consequence of the original Grothendieck's intension  to obtain that the topological space of the set of objects of the (for us derived) category $\textbf{C}$ is a \emph{scheme} (the topological space of open sets together with a commutative ring structure) in order to enlarge the notion of \emph{algebraic variety}, while in the Metacategory theory, representing \emph{all} categories, the unique common  space of each category  is composed by the simplest discrete collection of its objects (points) and oriented paths between them (the morphisms of the category).
Consequently, Grothendieck topology approach is not a generally valid model of the Metacategory's discrete topological space (that is, the abstract space of \emph{all} categories)\footnote{Moreover, "the introduction of categories as a part of the language of mathematics has made possible a fundamental, intrinsically categorical technique: the element-free definition of mathematical properties by means of commutative diagrams, limits and adjoints", (p. 3 in \cite{BaWe85}). So, for our scope of using only (co)limits and adjunctions, we do not need the "set-based" definition of  elements (i.e., the "points") as defined above by the special morphisms from terminal object \textbf{1}. One
of the main benefits of category theory is you do not have to do things in terms of elements unless it is advantageous to.}.
So, in what follows, we will introduce a more general model of the Metacategory space, first time used in \cite{Majk23s} but without a detailed explanation.
\section{A General Definition of Geometric Metacategory Space\label{sec:metacat}}
Let us consider the following  text taken from  \cite{DeBr11} (p. 1):

\emph{Albert Einstein once said, about "the world of our
sense experiences", "the fact that it is comprehensible
is a miracle" } (\cite{Eins36}, p. 351).\emph{ A few decades
later, another physicist, Eugene Wigner, wondered about
the unreasonable effectiveness of mathematics in the
natural sciences, concluding that "the miracle of the appropriateness
of the language of mathematics for the formulation
of the laws of physics is a wonderful gift which we
neither understand nor deserve"} (\cite{Wign60}, p. 14). \\\emph{At least three
factors are involved in Einstein’s and Wigner’s miracles: the
physical world, mathematics, and human cognition. }\\

In my recent work \cite{Majk23s} I considered not only the topological properties of category theory, but also its  geometry based in noncontinuous (discrete)  space composed by the  set of points (corresponding to the positions, in this  space, of the objects of a category) and oriented paths in such categorial space used to represent the morphisms of the category. So, in such a 3D categorial space, each object has a distinct position (point of this space) and the spatial relationship between them is defined by the oriented paths between object's positions (representing spatially the morphisms between the objects).

Thus, differently from the previous consideration of the "points" in the category theory (the Metacategory),  here, in an analogy with the approach used in \cite{DeBr11}, we consider the category theory by its characteristic commutative-diagrams \emph{visual graphical language} and hence we \emph{do not} consider the points as the elements of the object of a given category (which can not be applied to all categories; not all of them are well-pointed).

A discrete space is a particularly simple example of a  structure, one in which the points form a discontinuous sequence, meaning they are isolated from each other in a certain sense. For example, for any category the collection of its objects is a particular discrete subspace of such a category.

Let us provide a  3D representation of such  categorial space and in which way these points (corresponding to the objects of a category) are isolated, from the fact that we can consider the morphisms between the objects as the oriented paths (composed by the points $(x,y,z)$ in the 3D continuous space $Z$ as well).
So, let us consider the discrete object's space (O-space) of a given category $\textbf{C}$, as a collection of distinct points in the 2-dimensional $xy$   plain\footnote{This plain is introduced as a projective plain in which we usually visually represent the oriented graphs of categorial (commutative) diagrams.} $Z_0$ (with $z =0$), by introducing the mapping
\begin{equation} \label{eq:ObPos}
P:Ob_\textbf{C} \rightarrow Z_0
\end{equation}
such that for any two different objects $a,b \in Ob_\textbf{C}$ with $(x_1,y_1,0) = P(a)$ and $(x_2,y_2,0) = P(b)$ it holds that $(x_1-x_2)^2 +(y_1-y_2)^2 > 0$. And we can denote this discrete set of points in $Z_0$ plain as $Im(P)$ (image of the mapping $P$), and the continuous 3D subspace without this plain $Z_0$, $$Z\backslash Z_0 = \{(x,y,z)\mid x,y,z \in\mathbb{ R}, ~with~ z\neq 0\} $$
Consequently, the morphism's space (A-space) can be represented by the  3D subspace
$$A-space ~~\equiv ~~~(Z\backslash Z_0) \bigcup Im(P)$$
so that all points of each morphism have only the initial and final points of this morphism in the 2D plain $Z_0$ (which are the positions of the objects in $Im(P) \subset Z_0$), while the rest of the points of this  morphism, considered as an oriented path, are inside the 3D subspace $Z\backslash Z_0$ and hence are not in the plain $Z_0$. Thus, we consider that each path (without its initial and final positions) is totally above or bellow the plain $Z_0$, so that do not exist their intersections with the plain $Z_0$ which would introduce another points in it different from the points of the objects in $Im(P) \subset Z_0$.

 So, in the space of objects $Im(P)$ we have  the discrete space for a given category (objects are isolated) while in the  categorial 3D A-space the objects are not isolated but connected by the oriented 3D paths corresponding to the morphisms of a category.  In this way we represented the fact that, in the 2D discrete space $Im(P)$, the objects are "isolated from each other in a certain way".

Now, from the fact that the objects of the category can be eliminated and represented by their identity morphisms, the category can be represented geometrically by the 3D  A-space with the oriented paths inside it (representing the morphisms of the category) while the initial and final points of these paths inside the 2D discrete space $Im(P)$ represents the objects and not the morphisms. So,  in this real 3D space representation of the Metacategory space of the morphisms (paths) of a category, for the paths we have the space continuum and topologically we can describe this space (in the work of Ehresmann, Brouwer and Weil's \emph{Das Continuum})  by the class of its open sets (set of the points in the path without its initial and final position in the flat surface $Z_0$ that represent the objects).

In this way we can represent visually the \emph{diagrams} of a category (characteristic visual language representation in category theory) as a projection of the 3D A-space on the 2D plain $Z_0$ (both with the object's positions in $Z_0$).
That is,  we are able to define the categorial bidimensional \emph{commutative diagrams} (without unnecessary identity morphisms) obtained by this 3D A-space:
 \begin{definition} \label{def:comdiag} \textsc{Categorial Commutative Diagrams in A-spaces}:
  The commutative diagrams used as a graphic representation in category theory are the projection of the A-space into the plain surface $Z_0$.\\
  The representation of the categorial morphisms in commutative diagrams by simple bidimensional \textsf{linear vectors} can be obtained by the constraint that each path between two objects lies in the plain orthogonal to $Z_0$ passing these two source and target objects of such a oriented path.
  \end{definition}
Indeed, by this projection we obtain visible the morphisms in this plain surface $Z_0$ between the objects of a category exactly as we are using the commutative diagrams in the category theory.  By using the constraint for representation of categorial morphisms by \emph{linear vectors}  we obtain the relationship between the metacategory and geometry and vector spaces, as, for example, in diagrams in previous section. The equations representing the commutative diagrams are traduced into two \emph{equal vectors}.

With this representation of 3D metacategory space with defined constraints, we can now proceed to the process of "geometrization" of the category theory by definition of a kind of particular geometric algebra in next section,
obviously we have to consider the most known 3D geometric algebra of Clifford, and to examine their similar and also different properties.

For example, our 3D metacategory space of vectors is different from 3D exterior and Clifford algebras because we can not make the sum of any two vectors in this categorial geometric algebra because of partial composition of the morphisms in a category. Thus in what follows we will make a short introduction  to Grossmann exterior algebra and geometric (Clifford) algebra, which will be used in next section of definition of Cat-arrows space of vectors and its properties.
In mathematics, a
geometric algebra (also known as a real Clifford algebra which is an extension of elementary algebra to work with geometrical objects such as vectors) is built out of two fundamental operations, addition and the geometric product. Multiplication of vectors results in higher-dimensional objects called multivectors. The geometric product was first briefly mentioned by Hermann Grassmann. Clifford defined the Clifford algebra and its product as a unification of the Grassmann (exterior) algebra and Hamilton's quaternion algebra.
In the Grasmann algebra, the scalars and vectors have their usual interpretation, and make up distinct subspaces of a geometric algebra. Bivectors provide a more natural representation of the pseudovector quantities in vector algebra such as oriented area, oriented angle of rotation, torque, angular momentum and the electromagnetic field. A trivector can represent an oriented volume, and so on.
 An exterior algebra of a vector space $V$ is a graded associative algebra,
 \begin{equation}\label{extAlg}
     \bigwedge(V) =  \bigwedge^0(V) \bigoplus \bigwedge^1(V)\bigoplus \bigwedge^2(V)\bigoplus \bigwedge^3(V)\bigoplus ...
 \end{equation}
where $\bigoplus$ is direct sum, $ \bigwedge^0(V)$ the set of scalars, $\bigwedge^1(V) = V$, and $\bigwedge^2(V)$ the set of bivectors $v_1\wedge v_2$ for each two vectors $v_1,v_2 \in V$, where $\wedge$ is a noncommutative outer (wedge) product. Every element of $\bigwedge(V)$ is direct sum of $v_1\wedge ...\wedge v_k$, where $v_i \in V$ are vectors.
An element of this form is called k-blade (k-vector) and corresponds to the oriented parallelotope spanned by them with k-volume.

In linear geometry the decomposable k-vectors have geometric interpretations: the bivector $ u\wedge v$ represents the plane spanned by the vectors, "weighted" with a number, given by the area of the oriented parallelogram with sides $u$ and $v$. Analogously, the 3-vector $u\wedge v\wedge w$ represents the spanned 3-space weighted by the volume of the oriented parallelepiped with edges $u,v$, and $w$.

Geometric algebras are closely related to exterior algebras. The conventional definition of geometric algebra is carried out in the context
of vector spaces endowed with an inner product, or more generally a quadratic form.
In geometric  algebras we have the \emph{geometric product} $uv$ of two vectors that extends the exterior (outer) product by inner (scalar) vector product component, that is,
 \begin{equation}\label{extAlg1}
     uv =  u\cdot v + u\wedge v
 \end{equation}
 where $u\cdot v \in \bigwedge^0(V)$ is an inner product (ordinary scalar product) of vectors.

We consider here a vector space $V$ of arbitrary dimension over some field $K$. The Clifford algebra is a unital associative algebra that contains and is generated by a vector space $V$ over a field $K$, where $V$ is equipped with a quadratic form $Q : V \rightarrow K$. The Clifford algebra $Cl(V, Q)$ is the "freest" (most general) unital associative algebra generated by $V$subject to the condition
$ v^2 = Q(v)1$, for all $v\in V$,
where the product on the left is that of the algebra, and the $1$ is its multiplicative identity. Clifford algebra $Cl(V, Q)$  for $Q = 0$ reduces to exterior algebra.
If the vector space $V$ is real 3D space with orthonormal basis $\{\textbf{e}_1,\textbf{e}_2,\textbf{e}_3\}$ , then the Clifford (geometric) product (\ref{extAlg1}) yields the following relations for $1\leq i,j \leq 3$:
\begin{equation}\label{extAlg2}
    \textbf{e}_i^2 = 1,  ~~~~~~and ~~~\textbf{e}_i\textbf{e}_j = -\textbf{e}_j\textbf{e}_i ~~~~~~if~~ i \neq j
 \end{equation}
In what follows, in next sections,  we will show that the Cat-algebra based on an abstract Cat-arrows space\footnote{This space is composed by vectors but is not a vector space as in exterior and Clifford algebras.} $V$ of \emph{each} particular category (\emph{categorial space}), has similar "geometrical" properties. Such abstract categorial space is composed by discrete set of points (corresponding to the positions, in this abstract space, of the object of a given category) and vectors in this space are just the arrows of a category representing the oriented paths in such a space. In this way we can represent the \emph{diagrams} of a category (characteristic visual representation in category theory) as oriented graphs.

Note that the Clifford properties of his 3D geometric algebra are equally valid for our 3D A-space and  Definition \ref{def:comdiag}.
We will show that the vectors used for visual representation of the categorial morphisms do not generate a vector space used in exterior and Clifford algebras but a less general Cat-arrows space (categorial space of vectors). However, thanks to introduction of the norms for the morphisms, for which we can give different semantic representations,  we can have also the scalar multiplication of vectors, where the scalar modify the value of this weight for each concrete morphism, depending on any physical meaning (they are abstract concepts) of morphisms: considered as simple paths to cross in this 3D space and the length, necessary time to cross them, or the costs to cross them, etc.., if for example, the morphisms of category are the computation processes, functions, etc..

The fact that our Cat-arrow space of vectors defined in what follows, is only a very particular case of a more general Clifford algebra mainly is the result of the facts that not all morphisms in a given category can be composed, so also in their vector-representation we have only partial sum of these vectors (not generally we can not make the sum  any two two vectors derived from category morphisms).

I used such a categorial abstract space \cite{Majk24m} in order to show the global categorial symmetries and invariance in an analogy with the symmetries and invariance used in \cite{Majk23s} for the Einstein GR real geometry.
 In Einstein's conception of \emph{physical} spacetime (in a fixed time-instance) is not some pre-existing void in which matter, gravity and other forces of nature exist. It is real, substantial entity (not void) which is the gravitational tensor \emph{field} ($4\times 4$ metrics $\textbf{g}_{jk}$ identified as the gravitational potential) of matter sources. So, in any fixed time-instance, the curved 3D open space is the gravitational field.

 Symmetries play a fundamental role in physics because they are
related to conservation laws. This is stated in Noether's theorem
which says that invariance of the action under a symmetry
transformation implies the existence of a conserved quantity. Albert Einstein once said, about "\emph{the world of our
sense experiences}", and "\emph{the fact that it is comprehensible is a miracle}" (1936, p. 351). A few decades later, another physicist, Eugene Wigner, wondered about the unreasonable effectiveness of mathematics in the
natural sciences, concluding that "\emph{the miracle of the appropriateness of the language of mathematics for the formulation of the laws of physics is a wonderful gift which we neither understand nor deserve}" (1960, p. 14). At least three
factors are involved in Einstein’s and Wigner’s miracles\footnote{It is interesting fact that both of them considered the standard quantum mechanics (QM) as an \emph{incomplete} and statistical theory, so that for the description of the individual particles we needed a new more sophisticated theory with hidden variables. This was a starting point for me to develop this new and complementary part of QM from 2010 to 2019 published in three volumes \cite{Majk17,Majk17b,Majk17q} for completion of QM and unification with Einstein's General Relativity.}: the
physical world, mathematics, and human cognition.

We assume that in the absence of the \emph{adjunctions that generate the (co)limits} in a given category $\textbf{C}$ (in \cite{Majk23s} such adjunctions are considered as a categorial analog to the gravitational field) this A-space is a flat Euclidean subspace, and that if we have such adjunctions in this category that this $3D$ A-space is a curved GR space generated by the gravitation field represented by such categorial adjunctions. In this case, this 3D space curvature changes the point positions of the objects, that is, changes the mapping $P:Ob_\textbf{C} \rightarrow Z_0$.

 Consequently, in such an analogy with the real physical 3D space, the 3D categorial A-space can be  curved by the "categorial field" which is assumed in \cite{Majk23s} to be represented by the adjunctions. So in the presence of adjunctions that generate, for example, the (co)limits, categorial A-space (in which the time is frozen) is curved by these adjunctions, while in the absence of such adjunctions in a given category we can consider that A-space is a flat Euclidean 3-D space.

 In effect, this paper is only a more detailed explanation of the analogy of the invariances and symmetries in category theory introduced in \cite{Majk23s} and analog phenomena in Einstein cosmology.
This mathematical intuition, based on the "\emph{adjunction-as-filed}" paradigm, is summarized by table (2.53)  in Section 2.5.1of the book \cite{Majk23s}:
\begin{center}
 \begin{tabular}{|c|l|}
   \hline
    \textbf{Physics} Theory  & \textbf{Category} Theory\\
     \hline
 &\\
 Continuous time-space & Category $\textbf{C}$ as a 3D space with objects $p_i$ as points \\
with points $p_i$ & and arrows as directed paths between them\\
\hline
&\\
    Scalar field $\Psi$    & Adjunction $(\vartriangle, G,\varepsilon,\eta)$ with diagonal functor $\vartriangle:\textbf{C}\rightarrow \textbf{C}^{\textbf{J}}$\\
    of a given type $\textbf{J}$ &  of type defined by small index category $\textbf{J}$\\
\hline
&\\
    Lagrangian density $\L(\Psi,p_i)$    & Universal arrow $(\vartriangle(p_i),\eta(p_i))$ at a point $p_i$ \\
\hline
&\\
    Action $S =\int \Omega_g\L(\Psi,p_i)$    & Set of displacements of each point $p_i$ (path-integrals)\\
    &  represented by arrows with domain objects at $p_i$\\
\hline
&\\
    Principle of \emph{minimal} action    & Property that for each  cone of a diagram $d'$ in $\textbf{C}^{\textbf{J}}$\\
    &  there is a \emph{unique} arrow to the limit cone\\
\hline
&\\
    Euler-Lagrange equation    & Commutative diagram in  $\textbf{C}$ of the adjunction for\\
     with position variable $p_i$&   a fixed diagram $d'$ ("material object") in $\textbf{C}^{\textbf{J}}$ \\
     of a material object & and for each position variable $p_i$\\
 \hline
\end{tabular}
%
\end{center}
and the work of geometrization\footnote{We mean the displacement of the space points in the $Z_0$ plain of the A-space, caused by the curvature generated by the adjunctions (categorial fields) w.r.t. the flat space in absence of the adjunctions.} of the Category Theory (i.e., of the Metacategory), provided in this paper, gives more light to such intuition.

 It will be explained and analyzed in details in next section by definition of the abstract \emph{Cat-arrows space}, based on this real 3D Metacategory space representation  provided by the A-space above. This consideration provides the analogy between  categorial space and real Eintein's GR space in a fixed time instance.
 With such a kind of geometrization and the analogy between Physics and Metacategory space, it has been possible to investigate the symmetries and invariance in the Category Theory \cite{Majk23s} as well.

In this \emph{analogy} of the approach to the  relationship between Physics and Category theories, provided in in \cite{Majk23s}, we did not consider  the active transformations (as, for example, the Lorentz boosts). In the place of them, we  considered, in Section 2.5.1 in \cite{Majk23s}, another transformation from one to another reference frame with the constant velocity  $\overrightarrow{\textbf{v}}$ between them ("Galilean boosts"). Though the transformations are named for Galileo, it is absolute time and space as conceived by Isaac Newton that provides their domain of definition. \\
In next sections will be provided much more general theory for categorial space of vectors than in previous preprint, as invoked for future work in previous preprint \cite{Majk26a}.

\section{Categorial Space of Vectors Based on Oriented Paths}
In one of the first papers in topology, Leonhard Euler demonstrated that it was impossible to find a route through the town of K\"{o}nigsberg (now Kaliningrad) that would cross each of its seven bridges exactly once. This result did not depend on the lengths of the bridges or on their distance from one another, but only on connectivity properties: which bridges connect to which islands or riverbanks. This Seven Bridges of K\"{o}nigsberg problem led to the branch of mathematics known as graph theory.

"\emph{In accepting space as its object of investigation, geometry began to study relational structures instead of single figures or magnitudes (like triangles or conic sections). In this sense, the entire
structuralist approach of modern mathematics is grounded in this important shift of
perspective of eighteenth-century geometry, which (to use Cassirer’s words) turned
a classical geometry of substances (i.e. figures) into a geometry of functions
(structures)....\\
While the primary aim of classical geometry was the calculation of lengths,
areas and volumes of given figures, in a geometry of space the notions of position,
incidence or direction may play a central role.}" p. 2 of Introduction in \cite{DeRi15}

 The graph theory is fundamental part of the category theory, that is, each graph can be extended into a category as, for example, the small index categories $\textbf{J}$ derived from special graphs and used for co(limits) and categories of functors $\textbf{C}^{\textbf{J}}$.

Based on the axiomatic definition of the Metacategory in previous Section, in order to define the geometric properties of the category theory based on its  categorial A-space representation in previous section, we need to define what can be the
abstract \emph{space of vectors} of this geometry. The first point is to assume that each oriented path in such an abstract space (represented by an morphism of a category) is mathematically representable by a \emph{linear vector}, as specified by Definition \ref{def:comdiag}, in this categorial abstract space, and we will denote by $V$ the set of all vectors of a considered category: we will denominate it as the \emph{Cat-arrows space} to differentiate it from a common mathematical vector space.
\begin{definition} \label{def:vectors}
The categorial abstract space (Cat-arrows space) $V$ of any given category $\textbf{C}$ can be mapped into  real 3D space with orthonormal basis $\{\textbf{e}_1,\textbf{e}_2,\textbf{e}_3\}$  (used for 3D Clifford algebra as well) where $\{\textbf{e}_1,\textbf{e}_2\}$ are base vectors of the flat surface $Z_0$ of categorial commutative diagrams in Definition \ref{def:comdiag}:
\begin{enumerate}
  \item Each morphism $f\in Mor_\textbf{C}$, different from identity morphisms in this category\footnote{The identity morphisms can be omitted from the fact that they corresponds to the objects of the category and can be omitted from composition of morphisms in a category.}, is represented by a vector $f = a\textbf{e}_1+ b\textbf{e}_2$ with real constants $a,b \in \mathbb{R}$, and with its length $\|f\|$ equal to \emph{the time necessary to cross its path} by recommended  speed $v_0$ represented by real number 1.0
      in the 3D A-space\footnote{Note that $\|f\|$ is generally different from length $\sqrt{a^2+b^2}$ of the projection of this oriented path into plain surface $Z_0$. Another  example: if a morphism $f$ is a function or any process, as traveling for example, the $\|f\|$  can represent the average time for its execution.} .
  \item We define an auxiliary \emph{zero vector} in $V$, denoted by $\textbf{0} = \|\textbf{0}\|\textbf{e}_3$ in the space over the plain $Z_o$ (so that ist domain and codomain objects does not belong in plain $Z_0$), as vector representing the identity morphism of an object not in $Ob_\textbf{C}\in Im(P) \subset Z_0$, as specified by (\ref{eq:ObPos}), that is, with $dom(\textbf{0}) =cod(\textbf{0}) \notin Z_0$.

  \item The commutative and associative addition $\oplus$ of two given vectors $f$  and $g$ in $V$ corresponds to the partial operation of composition (constrained by the point 2 of the metacategory in  Section \ref{sec:SymMet}), defined as follows by considering that $\|\textbf{0}\| = 0$:
   \begin{equation} \label{eq:Comp-vectors}
         f\oplus g \equiv
 \left\{
    \begin{array}{ll}
        g, & \hbox{if $~~~f = \textbf{0} $}\\
     f, & \hbox{if $~~~g = \textbf{0} $}\\
    vector~ of~ morphism ~~ g\circ f \in Mor_\textbf{C}, & \hbox{if $~~~cod(f) = dom(g) $}\\
     vector ~of ~morphism ~~f\circ g\in Mor_\textbf{C}, & \hbox{if $~~~cod(g) = dom(f) $}\\
      non ~defined, & \hbox{~otherwise}
       \end{array}
  \right.
  \end{equation}
  so that if $g\circ f \in Mor_\textbf{C}$, and their vectors in plain $Z_0$ are $f = a_1\textbf{e}_1+ b_1\textbf{e}_2$ and $g = a_2\textbf{e}_1+ b_2\textbf{e}_2$, for reals $a_1,b_1,a_2$ and $b_2$, then the vector representation  $f\oplus g$ in $Z_0$ of composed morphisms in $\textbf{C}$ is equal to $f\oplus g = (a_1+b_1)\textbf{e}_1+ (a_2+b_2)\textbf{e}_2$ as standard sum of vectors $f$ and $g$.
  However, the semantics of operator $\oplus$ corresponds to the composition of the paths in A-space (the morphisms in category $\textbf{C}$), which projection into plain surface $Z_0$ define the vectors $f$ and $g$, and hence,
  \begin{equation} \label{eq:Comp-vectors2}
  \|f\oplus g \|\equiv \|f\| + \|g\|  \in \mathbb{R}
  \end{equation}
   This holds, from (\ref{eq:Comp-vectors}), also some of two vectors are zero vector $\textbf{0}$ because $\|\textbf{0}\| = 0$.\\
 If in some algebraic expression composed by more than one addition of vectors we obtain that anyone of them can not produce a vector in $V$, then all such vector's addition can not produce a well defined vector in $V$.
 \item The product of a vector $f$ with some positive scalar $\alpha$, is denoted by $\alpha f$ so that its norm is $\|\alpha f\| = \alpha\|f\|$.
      \end{enumerate}
  \end{definition}
 So, it is satisfied the property of zero vector, from (\ref{eq:Comp-vectors}), $\textbf{0}\oplus \textbf{0} = \textbf{0}$, from the fact that any  identity morphism  is idempotent w.r.t morphism compositions in a category.
  \\\\
 \textbf{Semantic Enrichment of Categories}:\\
 Notice that from (\ref{eq:Comp-vectors2}), the norm is not constrained to represent the simple geometric length of vectors in the plane $Z_0$ (in visual representation of commutative diagrams for which we are not interested for \emph{geometric length} of the vectors composing these visual diagrams) has the \emph{cumulative properties} over composition of morphisms in a category $\textbf{C}$ which define the paths in the 3D A-space. So, by this we introduce a new semantic enrichments for the morphisms of a category, by considering the composition of morphisms as definition of paths that in real application (execution of morphisms like  functions or processes that consume the real computational time, costs, or other cumulative physical phenomena, this additional semantics of morphisms, not definable in category, now is provided in category representation by these Cat-arrow space $V$ of vectors.\\ Let us consider, for example, a particular interpretation (semantics):
 \begin{example}
 Let us consider the seven bridges of K\"{o}nigsberg problem where each bridge is an object of category $\textbf{C}$ and the paths between them are represented by the morphisms, such that  starting from initial bridge we make the path to the last bridge by composition of the morphism (Euler demonstrated that at least for one bridge we have by this path pass two times). And we are interested for the total time to cross all these bridges by fixed recommended speed $v_0$ represented by real number 1.0. \\
 If this case for each morphism $f$ representing a particular path between two bridges (the objects in this category), its norm $\|f\| = 1.0\|f\|$ is just the time to cross the path of this morphism. And if the $n$ is the number of different morphisms that compose this total path $f = f_n\circ ...\circ f_2\circ f_1$, from first to last bridge, then the total time to cross this complete path is cumulative result, from (\ref{eq:Comp-vectors2}), equal to $\|f_n\otimes...\otimes f_2\otimes f_1\| = \sum_{i=1}^n \|f_i\|$.\\
 However,  if each path of morphisms $f_i$ is crossed by different average speed $v_i$, then we define the scalar $\alpha_i = \frac{v_0}{v_i} = \frac{1}{v_i}$, such that for vector $\alpha_i f_i$ we obtain that its norm is $\|\alpha_i f_i\| = \alpha_i\|f_i\| =\frac{1}{v_i} \|f_i\|$ is just the particular time to cross this path with its average speed.   So, for $\alpha < 1.0$ we travel with greater speed and for $\alpha > 1.0$ with smaller speed w.r.t. the recommended speed. In this case the norm of the linear vector in $V$, $g =\alpha_nf_n\otimes...\otimes \alpha_2f_2\otimes \alpha_1f_1$ will be real cumulative time to  pass from the first to last bridge, equal to
 $\|g\| = \sum_{i=1}^n \alpha_i\|f_i\|$.\\
 This demonstrate the importance of product of a vector with positive scalars, as defined in point 4 of Definition \ref{def:vectors}.\\
 In an analog way, we can consider the total cost or total energy needed to cross this complete path, and "crossing a morphism" can be interpreted by  consumption of time, money, or energy if this morphism is some function to be calculated or some process to be executed, for example. \\
 Consequently, we see how the "vectorization of categories" by definition the Cat-arrows 3D space of vectors, extends the categorial semantics with new external to category information representing by "norm of vector". This process is a kind of \emph{semantic enrichment} in Metacategory theory.
  \\$\square$
 \end{example}
 We can see that differently from the common (linear) vector spaces, in Cat-arrows space the addition of the vectors is not total but a partial operation, which is not defined for each pair of vectors, and for any vector in in Cat-arrows space $V$ we have no its inverse vector, and the scalars are restricted only on positive reals.

 In mathematics, a normed vector space or normed space is a vector space over the real or complex numbers on which a norm is defined. A norm is a generalization of the intuitive notion of "length" (in our case it is the time necessary to cross this oriented path) in the physical world. If $V$ is a vector space over $K$, where $K$ is a field of reals equal to $\mathbb {R}$  or to complex field
$\mathbb {C}$, then a norm on $V$ is a map $\|\_\|:V \rightarrow K$, satisfying the following four axioms:
\begin{enumerate}
  \item Non-negativity: for every $v\in V$, $\|v\|\geq 1$.
  \item Positive definiteness: for every $v\in V$, $\|v\|=0$ if and only if $v$ is zero vector.
  \item Absolute homogeneity: for every $v\in V$, and $\lambda \in K$, $\|\lambda v\| = |\lambda|\|v\|$.
  \item Triangle inequality: for every two vectors $v, u\in V$, $\|v\oplus u\| \leq \|v\|+ \|u\|$.
\end{enumerate}
 So, in our case case of Cat-arrows space with the partial sum operation of vectors $\oplus$ and the fact that we can not generally multiply the morphisms of a category by any kind of scalars, the Cat-arrows space is particular normed  space where the triangle inequality is just the equality in (\ref{eq:Comp-vectors2}) and the $K$ is field of \emph{positive} reals.

\section{Category Theory and Clifford 3D Algebra}
So, now we are able to define the inner (scalar) vector product of vectors in Cat-arrow space $V$, by considering that the angles between the vectors in the plain $Z_0$ are well defined:
\begin{definition} \label{def:scalarP}
For any two vectors $f$ and $g$ in $V$ (in Definition \ref{def:vectors}) with angle $\theta$ between them, their inner product '$\cdot$' can be defined in standard way by:
\begin{equation} \label{innerProduct}
 f\cdot g = (\|f\|\times \|g\|)\cos(\theta)
  \end{equation}
  where '$\times$' is the multiplication defined over real numbers semiring (it has all axioms of a ring excluding that of an additive inverse), $\mathbb{K} = (\mathbb{K}, 0,+,1,\times)$.
\end{definition}
Thus, from (\ref{innerProduct}), we obtain that $\textbf{0}\cdot \textbf{0} =\|\textbf{0}\|\times \|\textbf{0}\| = 0\times 0 =0$, and for every vector $f\in V$,
\begin{equation} \label{eq:V-baseZ}
f\cdot \textbf{0} = \textbf{0}\cdot f =0
\end{equation}
that is, the zero vector $\textbf{0} = \|\textbf{0}\|\textbf{e}_3$ is unique vector in $V$ orthogonal to \emph{all} other vectors in $V$.

We recall that \emph{nonassociative ring} is an algebraic structure that satisfies all of the ring axioms except the (full) associative property and the existence of a multiplicative identity. A notable example is a Lie algebra.
By considering the (partial) additive communicative operation $\oplus$ and (partial) commutative multiplicative operation $\cdot$ (inner product of vectors), we can represent the  algebraic structure $(V,\oplus,\cdot)$ of the Cat-arrows space $V$ by the following \emph{generalization} of the nonassociative ring (with a subset of axioms of the nonassociative ring):
\begin{definition} \label{def:V-ring}
The  algebraic structure $(V,\oplus,\cdot)$ of the Cat-arrows space $V$ is a commutative nonassociative ring,
So, we have only the following subset of axioms of the nonassociative ring:
\begin{enumerate}
  \item $(V,\oplus)$ is a partial (from (\ref{eq:Comp-vectors})) commutative monoid with identity element $\textbf{0}$ (called zero), such that:\\
      $(a\oplus b)\oplus c =a\oplus (b\oplus c)$, ~~if $~a\oplus b~$ and $~b\oplus c~$ are well defined vectors w.r.t (\ref{eq:Comp-vectors});\\
      $\textbf{0}\oplus a = a$;\\
      $a\oplus \textbf{0} =a$.
  \item Multiplication is commutative, such that:\\
  $\textbf{0}\cdot a =0$,\\
  $a\cdot \textbf{0} = 0$.
\end{enumerate}
\end{definition}
It is easy to verify that this definition is correct.

We can introduce the noncommutative partial substraction $\ominus$ of vectors in the Cat-arrows space, based on the noncommutative partial addition $\oplus$ in Definition \ref{def:vectors}, as follows:\\
For any two vectors $l$ and $g$ in $V$, if the morphism $f = l\circ g$ there exist in a category, that is there exist the vector $f = g \oplus l$ in $V$, then there exist also the distances

$d(f,g) = \|f\ominus g\| = \|l\|$ and $d(f,l) = \|f\ominus l\| = \|g\|$\\
which transforms this normed Cat-arrows space into a metric space.

%
%
So,  we are able to define also the outer (exterior or wedge) vector product of vectors in Cat-arrows space space $V$:
\begin{definition} \label{def:outerP}
The following applications of outer vector product in Cat-arrows space $V$ can be obtained:
\begin{itemize}
  \item For any two non-zero vectors $f$ and $g$ in $V$, obtained from the morphisms in category $\textbf{C}$  (in Definition \ref{def:vectors}), their outer product '$\wedge$' can be defined by bivector $f\wedge g$ , introduced by Grassmann exterior algebra, lying on the flat surface $Z_0$, so that:
\begin{equation} \label{outerProduct}
 f\wedge g = - g\wedge f
 \end{equation}
  \item For the zero vector $\textbf{0}$ in $V$, with $\|\textbf{0}\| = 0$, we obtain  that for each vector $f\in V$,

$\textbf{0} \wedge f = -f\wedge \textbf{0} = 0$,
so that also $\textbf{0}\wedge \textbf{0} =0$.
  \item If we use the base 3D vector $\textbf{e}_3$ which can not be the vector of any morphism in $\textbf{C}$, and is orthogonal to any vector obtained from the morphisms of category $\textbf{C}$, we can obtain also the bivectors not lying in the flat surface $Z_0$, and also non zero Grassmann trivectors (trivectors with volume different from zero).
\end{itemize}
\end{definition}
Note that  the bivectors $f\wedge g$ and $g\wedge f$ have the same area but with opposite direction (clockwise and counterclockwise) so that their sum is zero.

Consequently, we can define the graded Cat-algebra (Category algebra) as an extension of the Grasmann \emph{exterior} algebra, analogously to the 3D geometric (Clifford) algebra, over the Cat-arrows space $V$ with geometric vector product $fg = f\cdot g + f\wedge g$ introduced by (\ref{extAlg1}), defined over the semiring of real numbers $(\mathbb{R},0,+,1,\times)$.

Moreover, the geometric product of a non zero vector $f$ in $V$ by itself (we recall that $f$ is parallel to itself, so that $f\wedge f =0$), from (\ref{innerProduct}) and (\ref{outerProduct}), is equal to:
\begin{equation} \label{autoproduct}
f^2 = ff = f\cdot f = \|f\|^2 \geq 0
\end{equation}
So, from (\ref{eq:V-baseZ}), the only vector orthogonal to all other vectors
 is zero vector, and hence the inner product is nondegenerate.\\
 Heaving in mind this fundamental difference between the geometry of the Cat-algebra and the Clifford algebra, we have the following common consequences in these two geometric algebras:
\begin{coro} \label{Coro:c-algebra}
The Cat-algebra over the 3D Cat-arrows space $V$ with  geometric vector product over the semiring of real numbers $(\mathbb{R},0,+,1,\times)$, satisfies the Clifford algebra conditions:
\begin{enumerate}
  \item For each base vector $\textbf{e}_i \in V$ for $1\leq i\leq 3$, $\textbf{e}_i^2 = 1$, and $~\textbf{e}_i\textbf{e}_j = -\textbf{e}_j\textbf{e}_i$ if $~ i \neq j$.
  \item For any two mutually orthogonal vectors $f$ and $g$ in $V$, we obtain
  $fg = -gf$.
\end{enumerate}
\end{coro}
\textbf{Proof}: In fact,  the first point is the consequence of (\ref{extAlg2}).\\
For the second point, if $f$ and $g$ are orthodonal then $f\cdot g = g\cdot f = 0$, so that

$fg = f\cdot g + f\wedge g = f\wedge g = -g\wedge f = -g\wedge f- g\cdot f = -gf$.
\\$\square$
\section{Conclusion}
In section 2 we concluded that the Metacategory space, valid for all  categories, can not be defined by using well-know Grothendeick's approach with discrete ringed spaces.

In this way, the intuition expressed in \cite{Majk23s} that each category can be seen as a kind of abstract discrete space with points corresponding to its objects and oriented paths between the points of this abstract space corresponding to the non-identity morphisms of the category. In this paper we defined formally such abstract Cat-arrows space of the Metacategory (thus of each category) and, for the significant subset of the categories also the specific definition of their geometric Cat-algebra.

We recall that the Metacategory is a kind of associative non-commutative algebra of (partial) composition of the morphisms. Moreover, as a theoretical non set-based mathematical foundation, it is a general language in which all other algebras can be represented by particular categories.

Thus,  this "geometrization" of the Metacategory is a kind of techniques of applying geometrical constructions to algebraic problems, analog tho that used in the  branch of mathematics called \emph{algebraic geometry}, and hence we provide a relationship between the category theory (Metacategory) and this algebraic geometry based on Cat-arrows space of vectors. In this more general approach w.r.t. previous preprints, valid \emph{for all} categories, we finally  provided the general semantic definition of the norm $\|\_\|$ for the vectors in the 3D Cat-arrows space $V$, and how this definition enrich semantically the Metacategory theory, while this categorial geometric algebra is a particular restricted case of 3D Clifford algebras in which we can not have the sum for each pair of vectors in $V$ but only for that vectors that as morphisms can be composed in a given category and the product of scalars and vectors can be done only for positive scalars.

%

%

%
\newpage
\section{Appendix}
These are the figures: \\\\\\
 \begin{figure}
$\vspace*{-44mm}$
\centering{
 \includegraphics[scale= 1.0]{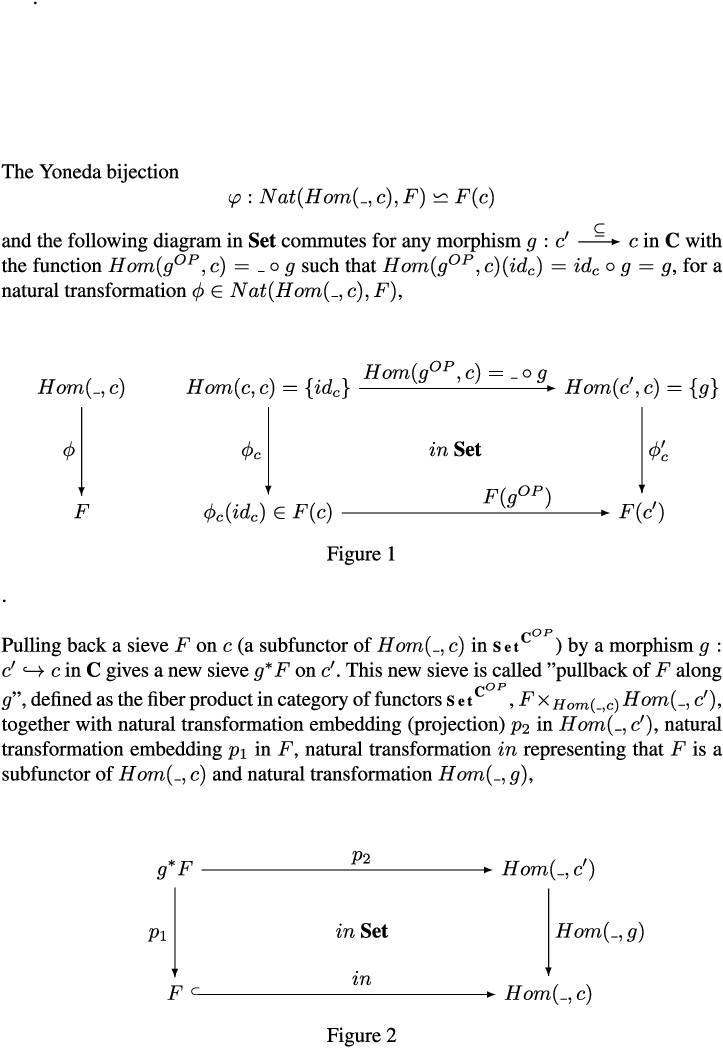}
 }
  $\vspace*{-11mm}$
 \end{figure}

\end{document}